\documentclass{article}
\usepackage[english]{babel}
\usepackage{amsmath,amsthm}
\usepackage{amsfonts}
\usepackage{authblk}
\usepackage{hyperref}
\usepackage{xcolor}

\usepackage{boldline}
\usepackage[margin=1in]{geometry}

\usepackage{graphicx}
\usepackage{chngpage}
\usepackage{array}
\usepackage{multirow}
\usepackage{slashbox}
\usepackage{makecell}

\usepackage{float}
\floatstyle{plaintop}
\restylefloat{table}

\usepackage[tableposition=top]{caption}
\usepackage{tabularx}
\usepackage{caption}
\usepackage{subcaption}

\usepackage{url}

\theoremstyle{definition}

\theoremstyle{remark}

\numberwithin{equation}{section}

\title{Marketing Mix Optimization with Practical Constraints}
\author{Hsin-Chan Huang\footnote{hsinchan.huang@nielsen.com, corresponding author} \hspace{0.25in}
Jiefeng Xu\footnote{jiefeng.xu@nielsen.com} \hspace{0.25in}
Alvin Lim\footnote{alvin.lim@nielsen.com}}
\affil{Precima, a NielsenIQ Company}
\date{\today}

\begin{document}
\maketitle

\begin{abstract}
In this paper, we address a variant of the marketing mix optimization (MMO) problem which is commonly encountered in many industries, e.g., retail and consumer packaged goods (CPG) industries. This problem requires the spend for each marketing activity, if adjusted, be changed by a non-negligible degree (minimum change) and also the total number of activities with spend change be limited (maximum number of changes). With these two additional practical requirements, the original resource allocation problem is formulated as a mixed integer nonlinear program (MINLP). Given the size of a realistic problem in the industrial setting, the state-of-the-art integer programming solvers may not be able to solve the problem to optimality in a straightforward way within a reasonable amount of time. Hence, we propose a systematic reformulation to ease the computational burden. Computational tests show significant improvements in the solution process.
\end{abstract}


\section{Introduction}
In today's business environment, a company may employ a mix of many marketing activities (e.g., free samples, discount coupons, weekly specials and advertisements) to attract new customers, enhance customer loyalty, and induce sales. Given a marketing plan with a mix of marketing activities, each activity may be allowed a spend level that results in corresponding outcome. The company often needs to adjust the resource allocated to marketing activities due to various reasons (e.g., marketing budget change, change in scope of targeted audience, geographies and products, product demand and/or supply changes, and activity cost and/or performance changes), aiming to optimize the revenue, profit, or other goals of the company. In practice, the company may desire such adjustments to be planned and executed by considering resource availability for executing changes, contractual obligations, etc. Due to minimum incremental investments required for certain marketing activities (e.g., minimum cost of a targeted rating point in TV advertisements or minimum cost of implementing a promotional discount), there is a need to impose a minimum change in spend (investment) constraint for a marketing activity. Since resources will be required to implement changes in marketing activities and resources available are limited, a constraint on the maximum to number of activities with spend changes is also imposed. Therefore, to guarantee that the outcome of a marketing mix optimization (MMO) is feasible for practical implementation, there is a need to solve the MMO problem with these added constraints.

Given an existing marketing plan (typically from a previous period), our MMO problem aims to adjust spend allocations on marketing activities with an objective of maximizing revenue subject to four sets of constraints: bounds on activity spend, total budget limits, the minimum change in spend for an activity, and the maximum number of activities associated with spend change. In the objective function of our problem, the revenue function is estimated as a (concave) separable quadratic function in terms of the activity investments, assuming the revenue response from investments in each activity is independent of one another. As a matter of fact, the marginal return on a marketing activity investment may not be constant, but rather diminishing. Such a plateau curve on the return of a marketing investment can be implemented using the left-half of a concave quadratic function together with an upper bound on the investment. In addition to the two constraints discussed earlier, when adjusting the spend of an activity from its baseline (previous period) level, the company may consider other business constraints. The activity spend may be restricted to a range to capture what is realistically implementable and to avoid a dramatic change from its baseline level. Also, the total marketing budget proposed may be limited based on the company's financial plan.

An intuitive way to formulate such a resource allocation problem is to set activity spend/investments as decision variables. By incorporating only constraints on activity spend bounds and total budget limits, the optimization problem becomes a quadratic programming (QP) with an objective function that maximizes the (concave) separable quadratic revenue function. It is well known that if the Q matrix of a QP is positive semi-definite, the problem can be solved in (weakly) polynomial time using the ellipsoid method~\cite{Kozlov_et_al_1980}. A commercial QP solver (e.g., Gurobi~\cite{Gurobi_2020}) can solve the problem quite efficiently. However, once we incorporate the minimum change in spend and maximum number of changes constraints mentioned above, we may need to introduce binary variables to take care of the discrete nature of the decision variables. Hence, the optimization problem becomes a mixed integer nonlinear program, which is well known to be NP-hard and generally cannot be guaranteed to solve in polynomial time.

In this study, we propose to employ a semi-continuous variable to represent a certain type of semi-continuity of an activity spend level and utilize its property to facilitate the solution process. A semi-continuous variable is defined as a variable whose value is either 0 or in an interval between a lower bound and an upper bound~\cite{Williams_2013, Sun_et_al_2013}, where 0 can be within the bound range. To align our activity spend level with the standard definition of the semi-continuous variable, we choose to use change in spend for an activity from the baseline level as the decision variable. The change in spend (adjustment) then becomes a semi-continuous variable. However, the minimum change requirement further splits the one continuous range into (at most) two disconnected continuous ranges. We present a systematic way to formulate the minimum change requirement by exploiting the semi-continuous nature of the variable. The constraint for the maximum number of changes in our problem is known as a cardinality constraint in the literature.

Optimization problems with semi-continuous variables and cardinality constraint have received much attention in the literature as there exists many applications in the real world, e.g., portfolio selection, and compressed sensing and subset selection~\cite{Sun_et_al_2013}. In particular, there exists quite a few papers discussing how to solve a variety of portfolio selection problems. Such problems may include minimizing the risk of a portfolio, represented by a positive semi-definite quadratic function, subject to a typical minimum gain constraint and several other linear constraints, resulting in a mixed integer quadratic program (MIQP). Our MMO problem with minimum change in spend and maximum number of changes shares a similar model structure with the portfolio selection problems, especially those with the minimum buy-in requirement. Inspired by the perspective reformulation Frangioni and Gentile~\cite{Frangioni_et_al_2006, Frangioni_et_al_2009} used to solve optimization problems with semi-continuous variables, where an MIQP is reformulated as a mixed integer second-order cone program, this paper proposes an alternative formulation for our MMO problem to ease the computational burden. The model formulation and research findings in this paper can be readily applied to problems with similar model structure and constraint settings.

The rest of the paper is organized as follows. Section~2 reviews literature on optimization models in marketing planning and optimization models with semi-continuous variables and cardinality constraint, with a focus on portfolio selection problems. Section~3 presents a systematic way to model semi-continuity on variables required to have a minimum change and proposes an alternative to model our MMO problem. Section~4 reports on computational tests comparing the solution outcomes of various formulations. Section~5 analyzes the computational results and discusses the pros and cons of the varying formulations. Finally, Section~6 concludes with highlights on the findings of this research.

\section{Literature Review}

Optimization methods and tools are valuable in providing significant performance improvement in the marketing efforts of today's business. Fischer et al.~\cite{Fischer_et_al_2011}, for example, constructed a dynamic optimization model to allocate marketing budget across countries, products and marketing activities, and proposed a heuristic to obtain a feasible solution that is close to optimal. They implemented it at Bayer, a large multi-national pharmaceutical and chemical firm, and estimated the profit improvement potential to be more than 50\% which is worth hundreds of millions of euros. In general, a decision support system for marketing consists of two stages: estimating the response of related variables and optimizing the performance over some of the variables. Such a system has been discussed extensively in the literature. For example, Balachandran and Gensch~\cite{Balachandran_et_al_1974} estimated sales as a function of the marketing mix and environmental variables and formulated the budget allocation problem as a geometric programming problem. More recently, Pferschy et al.~\cite{Pferschy_et_al_2014} estimated the direct and joint effects of different media in a promotional campaign and formulated the advertising media selection problem as a quadratic knapsack model. The quality of results from an optimization model is dependent on the response model. Albers~\cite{Albers_2012} argued that the output of the response model should be optimizable, and provided recommendations on optimizable and implementable aggregate response modeling for marketing decision support. For an overview of optimization methods in marketing, Mantrala and Kanuri~\cite{Mantrala_et_al_2018} conducted a survey over the last 70 years of research, where they summarized each example with the decision variables, the constraints, the sales response function and the objective function.

Our MMO problem, however, is unique to what exists in the marketing literature but is closely related to portfolio selection problem with the minimum buy-in threshold and limited total number of assets, which belongs to the family of mathematical programs with semi-continuous variables and cardinality constraint. Hence, we provide a literature review on the techniques used to solve mathematical programs with semi-continuous variables, especially in the context of perspective reformulation. Portfolio selection problems (especially a variety of Markowitz's mean-variance (MV) portfolio selection models) have received much attention in the field of finance and financial engineering. One popular mean-variance model is to minimize the volatility of portfolio returns (risk) subject to a guaranteed return, full capital invested and non-negative investment~\cite{Markowitz_1952}, the model is usually formulated as an MIQP. In addition to this basic model, many researchers have been adding more real-life constraints and proposing efficient ways to solve these more practical models. For example, Bienstock~\cite{Bienstock_1996} presented a computational experience with a branch-and-cut algorithm to solve the MV model with a cardinality constraint. Bertsimas and Shioda~\cite{Bertsimas_et_al_2009} described a branch-and-bound based algorithm with the use of Lemke's method to solve the MV model with the minimum buy-in threshold and cardinality constraint. Cesarone et al.~\cite{Cesarone_et_al_2015} studied the MV model with the minimum buy-in threshold and cardinality constraint, and two other variations. Instead of minimizing the volatility of return, one of their variations minimizes conditional value-at-risk and the other minimizes mean absolute deviation. They proposed an approach based on a reformulation as a standard QP and presented some theoretical results to solve the MV with minimum buy-in threshold and cardinality constraint.

One of the effective techniques to solve a mixed integer convex nonlinear programming problem with semi-continuous variables and a separable objective function is via the perspective relaxation proposed by Frangioni and Gentile~\cite{Frangioni_et_al_2006}. They constructed a convex envelope of the univariate function over a semi-continuous variable by using its perspective function. By exploiting the characterization of the perspective function, they derived a family of valid inequalities (referred to as perspective cuts, P/C) and incorporated them into a typical branch-and-bound algorithm, resulting in substantial performance improvement in solving a unit commitment problem. In addition to P/C reformulation, one can also reformulate the whole problem as a second-order cone program (SOCP) if the problem has a separable quadratic objective function. Frangioni and Gentile~\cite{Frangioni_et_al_2009} conducted a computational comparison of P/C versus SOCP reformulations on portfolio selection and unit commitment problems. They concluded that P/C reformulation seems to perform better due to much more efficient re-optimization capabilities of branch-and-bound algorithms compared to interior point algorithms at that time. Frangioni and Gentile~\cite{Frangioni_et_al_2007} extended the use of perspective relaxation to problems with non-separable quadratic objective functions. They proposed two methods to decompose the one positive semi-definite Q matrix into two positive semi-definite matrices, one generic and the other diagonal. One method is more efficient involving the calculation of the minimum eigenvalue and the other is more costly requiring the solution of a semi-definite program. Based on testing both methods using P/C reformulation on an MV model, they concluded that the more costly method significantly outperforms the other. Zheng et al.~\cite{Zheng_et_al_2014} later proposed to decompose the Q matrix by solving a new semi-definite program claiming that their decomposition gives the tightest continuous relaxation of the perspective reformulation of quadratic programming problems with minimum threshold and cardinality constraints.

There are many real world applications for the mixed integer nonlinear programming (MINLP) problems with semi-continuous variables in addition to portfolio selection and unit commitment problems. Gunluk and Linderoth~\cite{Gunluk_et_al_2012} discussed a variety of practical MINLP problems whose relaxation can be strengthened via the perspective reformulation, including separable quadratic uncapacitated facility location problem, network design problems with congestion constraints, scheduling problems with controllable processing times, unit commitment problems, stochastic service system design problems, and portfolio selection problems. Sun et al.~\cite{Sun_et_al_2013} surveyed developments in solving MINLP application problems with semi-continuous variables and cardinality constraint, including problems in production planning, portfolio selection, compressed sensing, and subset selection in regression.

While the research in MINLP applications is quite extensive, all the applications considered semi-continuous variables to be either 0 or within a continuous range (i.e., $x \in \{0\} \cup [\alpha, \beta]$). However, in our MMO problem, the change (either positive or negative) in an activity spend also requires to be at least no less than a certain threshold (i.e., $x \in \{0\} \cup [\alpha, \beta] : |x| \geq \delta$). Using an analogy to portfolio selection problem, we could allow either long selling or short selling a security with minimum buy-in threshold. We can put an absolute value operator to the semi-continuous variable to satisfy such a requirement. However, this will result in a nonlinear constraint to the semi-continuous variable. In this paper, we extend the perspective relaxation to the semi-continuous variable with the absolute value operator.

\section{Optimization Problem and Proposed Models}
To reiterate, our MMO problem is to maximize revenue subject to four groups of constraints: (a) bounds on activity spend, requiring the activity spend to be within a certain range; (b) total budget limit, requiring the total of all activity spends to be within a baseline budget; (c) minimum change in activity spend, desiring the change in activity spend (if any) to be at least a certain amount (see below); and, (d) maximum number of activities with spend change, limiting the required actions from a baseline plan. The decisions are which activities to change investments and by how much for each activity. In practice, there may be other linear constraints based on certain business requirements, which can be added to the model without suppressing the benefit of the following proposed reformulation. Without loss of generality, we consider only the four groups of constraints to keep the model succinct. In constructing the optimization model, instead of activity spend, we use the change in activity spend as the decision variable to facilitate later model reformulation. For the objective function, we assume that the revenue is a (concave) separable quadratic function of spend change. Furthermore, we assume that the coefficients of the second-order terms are non-positive so that the optimization problem is indeed a concave maximization problem.

By incorporating the minimum change in activity spend and the maximum number of activities with spend change constraints, we first formulate the MMO problem as an MINLP. We further present a systematic way to model the semi-continuity due to the minimum change requirement. With the (concave) separable quadratic revenue function and the pre-processing of the minimum change constraint and activity spend bounds, we can reformulate the problem as an MIQP. However, solving the MIQP to optimality directly using standard commercial solvers may be time-consuming and impractical. We then propose an alternative model of the optimization problem as a mixed integer second-order cone program (MISOCP), aiming to shrink the optimality gap and/or reduce the solution time. In the next several subsections, we first elaborate on the modeling of the semi-continuity and then explain the reformulation of the optimization problem in detail.

\subsection{The Minimum Change in Activity Spend}
When the spend on activity $i$, $y_i$, changes from its baseline level $s_i$, the activity spend range $[l_i, u_i]$ and minimum change requirement $\left |y_i - s_i \right |\geq \delta_i$ are imposed on the spend change $x_i$. We assume the baseline level is within the activity spend range and therefore 0 is feasible for the spend change.

\begin{figure}[h!]
  \includegraphics[width=\linewidth]{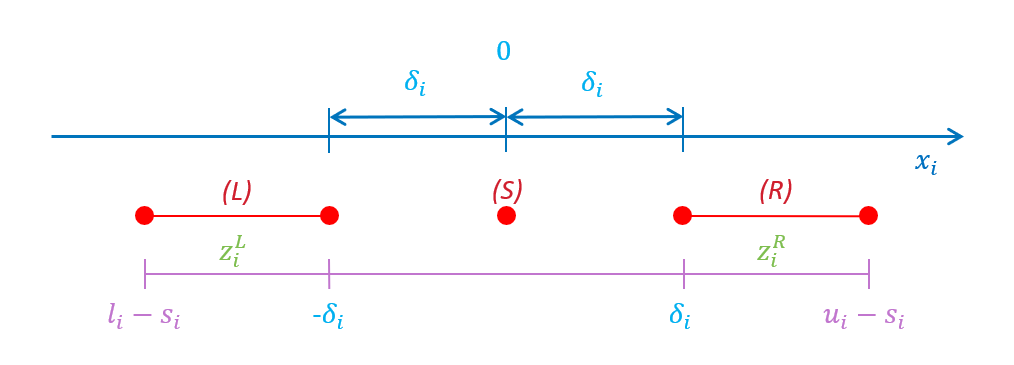}
  \caption{Illustration of the spend change of activity $i$, $x_i$, bounded by its activity spend baseline level $s_i$ and spend range $[l_i, u_i]$, the minimum change $\delta_i$, and the resultant feasible regions $L$, $S$ and $R$, shown by red segments.}
  \label{fig:MinimumChangeIllustration}
\end{figure}

The introduction of the minimum change constraint causes the feasible spend change $x_i$ to become discontinuous. Without such a requirement, the change can be any point within the range $[l_i - s_i, u_i - s_i]$. With the minimum change requirement, there are at most three possible regions that the spend change may fall into: region $L$, decrease in spend; region $S$, no change; and, region $R$, increase in spend, as shown in Figure~\ref{fig:MinimumChangeIllustration}. Note that regions $L$ and $R$ may or may not exist (feasible or infeasible), depending on the relative positions of the baseline level $s_i$ to the activity spend range $[l_i, u_i]$, and the minimum change ($\delta_i$).

Figure~\ref{fig:MinimumChangeIllustration} illustrates the scenario with the most (three) possible feasible regions for the spend change of activity $i$. The bounds of regions $L$ and $R$ may vary under different scenarios. Taking region $L$ as an example, if the activity spend lower bound minus the baseline level is no larger than the negative minimum change (i.e., $l_i - s_i \leq - \delta_i$), the lower bound of region $L$, $l_i^L$, is the activity spend lower bound minus the baseline level, and its upper bound $u_i^L$ is the negative minimum change (i.e, $[l_i^L, u_i^L]=[l_i - s_i, - \delta_i]$). On the other hand, if the activity spend lower bound minus the baseline level is larger than the negative minimum change, region $L$ does not exist. Similar relationships can be derived for region $R$. For region $S$, since we assume the baseline level is within the activity spend bound (i.e., $l_i \leq s_i \leq u_i$), the bound of region $S$ always exists and is 0. Based on these three regions, we introduce two binary variables ($z_i^L$, $z_i^R$, where the superscription denotes the region) to indicate which region that the optimal spend change falls into. Table~\ref{tab:MC_change} discloses all possible relationships between region bounds and binary variable values for the three regions under different scenarios.
\begin{table}[h!]
\footnotesize
\setlength{\extrarowheight}{0.1in}
\begin{adjustwidth}{-1in}{-1in} 
\begin{center}
    \begin{tabular}{| c | c | l | l | l |}
    \hline
    Region & Indicator Variable & Scenario & Variable Value & Region Range \\ \hline
    \multirow{2}{*}{$L$} & \multirow{2}{*}{$z_i^L$} & $l_i - s_i \le - \delta_i$ & $\in \{0, 1\}$ & $[l_i - s_i, - \delta_i]$ \\ \cline{3-5}
    &  & $ - \delta_i < l_i - s_i$ & $\in \{0\}$ & N/A \\ \hline
    \multirow{1}{*}{$S$} & \multirow{1}{*}{$1 - z_i^L - z_i^R$} & $l_i \le s_i \le u_i$ & $\in \{0, 1\}$ & $[0, 0]$ \\ \hline
    \multirow{2}{*}{$R$} & \multirow{2}{*}{$z_i^R$} & $\delta_i \le u_i - s_i$ & $\in \{0, 1\}$ & $[\delta_i, u_i - s_i]$ \\ \cline{3-5}
    &  & $u_i - s_i < \delta_i$ & $\in \{0\}$ & N/A \\ \hline
    \end{tabular}
\end{center}
\end{adjustwidth}
\caption{Possible regions for the spend change of activity $i$, $x_i$, under different scenarios, where $s_i$ is the activity spend baseline level, $l_i$ and $u_i$ are the lower and upper bounds of activity spend, $\delta_i$ is the minimum change, and $z_i^L$ and $z_i^R$ are indicator variables.}
\label{tab:MC_change}
\end{table}

We therefore formulate the activity spend bound and the minimum change in spend together as follows, which we will subsequently incorporate into the optimization model:
\begin{subequations} \label{MC_const}
\begin{align}
	& l_i^L z_i^L + l_i^R z_i^R \leq x_i \leq u_i^L z_i^L + u_i^R z_i^R, \quad \forall i \in I \label{MC_con1}\\
	& z_i^L + z_i^R \leq 1, \quad \forall i \label{MC_con2}\in I\\
	& z_i^L \in B_i^L, z_i^R \in B_i^R, \quad \forall i \in I \label{MC_con3}
\end{align}
\end{subequations}
In the above formulation, constraint (\ref{MC_con1}) specifies the lower bound and upper bound of the spend change. Constraint (\ref{MC_con2}) ensures the spend change falls into one and only one of the three regions. Moreover, constraint (\ref{MC_con3}) requires the indicator variables be binary. Note that we specify the possible values of the two binary variables in $B_i^L$ and $B_i^R$ according to Table~\ref{tab:MC_change}.

An alternative to the above formulation may use a constraint that restricts the absolute value of the spend change to be no less than the minimum change if the activity spend does not stay at the baseline level. However, it will result in an MINLP due to the use of absolute value operator, which we will discuss in the next subsection.

\subsection{The Convex Envelope of a Univariate Function Over a Semi-continuous Variable with Absolute Value Operator}
Frangioni and Gentile~\cite{Frangioni_et_al_2006} constructed a convex envelope of the univariate function over a semi-continuous variable and proposed to use it as a building block to reformulate their original problem involving semi-continuous variables, attempting to obtain a tighter continuous relaxation. They considered the feasible value of a semi-continuous variable as either 0 or within a continuous range. In this paper, we extend their approach to a semi-continuous variable with the absolute value operator, i.e., a variable is either 0 or within multiple continuous ranges.

Considering the minimization problem in~\cite{Sun_et_al_2013}, we add the absolute value operator to the semi-continuous variable in the disjunctive constraint and use a slightly different notation as follows:
\begin{subequations} \label{co2_1}
\begin{align}
	\quad \min_{x, t} \quad & f(x) + dt \label{co2_1_obj} \\
	\textrm{s.t.} \quad & \alpha t \leq | x | \leq \beta t \label{co2_1_con1} \\
	& t \in \{0,1\} \label{co2_1_con2}
\end{align}
\end{subequations}
In the model, $x$ is the semi-continuous variable whose value is either 0 or a point within $[\alpha, \beta]$ or $[-\beta, -\alpha]$, and $t$ is the indicator variable whose value is 1 if $x$ is nonzero and 0 otherwise. The objective function consists of two parts, one is a continuous function of $x$ and the other is a flat cost $d$ if $x$ is nonzero.

Without the absolute value operator, the perspective reformulation approach for the minimization problem proposed by Frangioni and Gentile~\cite{Frangioni_et_al_2006} is to replace the continuous part of the objective function with its perspective function (i.e., $f(x) \to t f(x/t)$). By adding the absolute value operator, the feasible value of the semi-continuous variable changes from $x \in \{ 0 \} \cup [l^R, u^R]$ (\{0\} plus a continuous range) to $x \in \{ 0 \} \cup [l^L, u^L] \cup [l^R, u^R]$ (\{0\} plus two continuous ranges), where $[l^L, u^L]$ and $[l^R, u^R]$ represent $[-\beta, -\alpha]$ and $[\alpha, \beta]$, respectively, for the rest of the subsection. To remove the absolute value operator, we introduce an additional binary variable $\overline{t}$ and rewrite the problem as follows:
\begin{subequations} \label{co3_1}
\begin{align}
	\quad \min_{x, t} \quad & f(x) + dt \label{co3_1_obj} \\
	\textrm{s.t.} \quad & \left[ l^L \left( 1 - \overline{t} \right) + l^R \overline{t} \right] t \leq x \leq \left[ u^L \left( 1 - \overline{t} \right) + u^R \overline{t} \right] t \label{co3_1_con1} \\
	& \overline{t} \leq t \label{co3_1_con2} \\
	& t, \overline{t} \in \{0,1\} \label{co3_1_con3}
\end{align}
\end{subequations}
In the above formulation, constraint~(\ref{co3_1_con2}) tightens the feasible region by excluding the point $\left( t, \overline{t} \right) = (0, 1)$, which makes $x = 0$ same as $\left( t, \overline{t} \right) = (0, 0)$ does. Given $t = 1$ in constraint~(\ref{co3_1_con1}), the feasible range of $x$ changes from $[l^L, u^L]$ to $[l^R, u^R]$ as $\overline{t}$ goes from 0 to 1. By specifying the value of $\overline{t}$, the two square brackets become specific lower and upper bounds for $x$ and the problem reduces to the minimization problem in~\cite{Sun_et_al_2013}, where its objective function can be rewritten using its perspective function to provide a tighter continuous relaxation. Note that, in general, the two continuous ranges need not be symmetric with respect to 0 (and can even overlap) although we derived it by removing the absolute value operator from the semi-continuous variable.

We now develop a tighter formulation for the nonlinear constraint~(\ref{co3_1_con1}) by employing a similar technique as in~\cite{Sun_et_al_2013}. Considering the lower bound of $x$ in the constraint, the feasible points $\left( t, \overline{t}, x \right)$ are $(0, 0, 0)$, $(1, 0, l^L)$, and $(1, 1, l^R)$. Any point in the smallest triangular area formed by these three points can be represented as:
\begin{align}
	\left( 1 - \theta \right) \left( 0, 0, 0 \right) + \theta \left( 1, \tilde{t}, \tilde{x} \right) = \left( \theta, \theta \tilde{t}, \theta \tilde{x} \right), \nonumber
\end{align}
where $\theta \in [0, 1]$, $\tilde{t} \in [0, 1]$, and $\tilde{x} = l^L \left( 1 - \tilde{t} \right) + l^R \tilde{t}$. Letting $t = \theta$ and $\overline{t} = \theta \tilde{t} = t \tilde{t}$, we have $\tilde{t} = \frac{\overline{t}}{t}$ and $0 \leq \frac{\overline{t}}{t} \leq 1$ (or $0 \leq \overline{t} \leq t$). Then, we have $\theta \tilde{x} = t \left[ l^L \left(1 - \frac{\overline{t}}{t} \right) + l^R \left( \frac{\overline{t}}{t} \right) \right] = l^L \left( t - \overline{t} \right) + l^R \overline{t}$. Hence, the lower bound of $x$ can be represented as:
\begin{align}
	& l^L (t - \overline{t}) + l^R \overline{t} \leq x \nonumber \\
	& t, \overline{t} \in [0, 1]. \nonumber
\end{align}
By applying a similar logic to the upper bound of $x$, problem (\ref{co3_1}) can then be equivalently rewritten as follows.
\begin{subequations} \label{co4_1}
\begin{align}
	\quad \min_{x, t} \quad & f(x) + dt \label{co4_1_obj} \\
	\textrm{s.t.} \quad & l^L \left( t - \overline{t} \right) + l^R \overline{t} \leq x \leq u^L \left( t - \overline{t} \right) + u^R \overline{t} \label{co4_1_con1} \\
	& \overline{t} \leq t \label{co4_1_con2} \\
	& t, \overline{t} \in \{ 0, 1 \} \label{co4_1_con3}
\end{align}
\end{subequations}

In other words, constraint (\ref{co3_1_con1}) is linearized by removing the multiplier $t$ outside the brackets and also replacing 1 with $t$ inside the brackets, which provides the convex hull of the nonlinear constraint. Moreover, we can perform variable transformation to replace $\left( t - \overline{t} \right)$ and $t$ with $z^L$ and $z^R$, and rewrite problem~(\ref{co4_1}) as follows.
\begin{subequations} \label{co5_1}
\begin{align}
	\quad \min_{x, t} \quad & f(x) + dt \label{co5_1_obj} \\
	\textrm{s.t.} \quad & l^L z^L + l^R z^R \leq x \leq u^L z^L + u^R z^R \label{co5_1_con1} \\
	& z^L + z^R \leq 1 \label{co5_1_con2} \\
	& z^L + z^R = t \label{co5_1_con3} \\
	& z^L, z^R \in \{ 0, 1 \} \label{co5_1_con4}
\end{align}
\end{subequations}
The above transformation enables the values of $z^L$ and $z^R$ directly indicate which region $x$ is in ($[l^L, u^L]$, $[l^R, u^R]$, or $\{0\}$), so we can use it to impose region-related constraints for our MMO problem. Then, we are able to replace the continuous part of the objective function with its perspective function to obtain a tighter continuous relaxation as follows since problems (\ref{co3_1}) and (\ref{co5_1}) are equivalent.
\begin{subequations} \label{co5_2}
\begin{align}
	\quad \min_{x, t} \quad & t f \left( \frac{x}{t} \right) + dt \label{co5_2_obj} \\
	\textrm{s.t.} \quad & l^L z^L + l^R z^R \leq x \leq u^L z^L + u^R z^R \label{co5_2_con1} \\
	& z^L + z^R \leq 1 \label{co5_2_con2} \\
	& z^L + z^R = t \label{co5_2_con3} \\
	& z^L, z^R \in \{ 0, 1 \} \label{co5_2_con4}
\end{align}
\end{subequations}
Furthermore, we extend the above derivation from \{0\} plus two continuous ranges to \{0\} plus any number of continuous ranges (say $n$ ranges) by replacing~(\ref{co5_2_con1}) to (\ref{co5_2_con4}) with the following constraints.
\begin{align}
	& l^1 z^1 + l^2 z^2 + \dots + l^n z^n \leq x \leq u^1 z^1 + u^2 z^2 + \dots + u^n z^n \nonumber \\
	& z^1 + z^2 + \dots + z^n \leq 1 \nonumber \\
	& z^1 + z^2 + \dots + z^n = t \nonumber \\
	& z^1, z^2, \dots, z^n \in \{ 0, 1 \} \nonumber
\end{align}
We thereby use formulation~(\ref{co5_2}) as a building block to formulate the MMO problem, aiming for a tighter continuous relaxation to accelerate the solution process.

\subsection{Optimization Models}
Before developing the marketing mix optimization model, we first define a few notations.  Let $I$ be the set of marketing activities and let $i \in I$ represent a marketing activity. The following input data and parameters are required to describe the optimization models:
\begin{adjustwidth}{0.25in}{0.25in} 
\begin{center}
\begin{tabular}{l c p{4.2in}}
	&&\\
	$f_i (x_i)$ & : & revenue of activity $i$ as a function of its spend change $x_i$; $f_i (x_i) = \theta_i x_i^2 + \phi_i x_i + \psi_i$\\
	$\theta_i$ & : & coefficient of the second-order term in the revenue function for activity $i$\\
	$\phi_i$ & : & coefficient of the first-order term in the revenue function for activity $i$\\
	$\psi_i$ & : & constant term in the revenue function for activity $i$\\
	
	$s_i$ & : & baseline spend of activity $i$\\
	
	$\rho$ & : & upper bound of total activity budget in percentage of total baseline activity budget\\
	
	$l_i$ & : & spend lower bound of activity $i$\\
	$l_i^L$ & : & region $L$ lower bound of activity $i$\\
	$l_i^R$ & : & region $R$ lower bound of activity $i$\\
	$u_i$ & : & spend upper bound of activity $i$\\
	$u_i^L$ & : & region $L$ upper bound of activity $i$\\
	$u_i^R$ & : & region $R$ upper bound of activity $i$\\
	
	$n$ & : & total number of activities (i.e., $n = |I|$)\\
	$m$ & : & maximum number of activities with spend change\\
	$\delta_i$ & : & minimum change in spend of activity $i$\\
	
	$M$ & : & a sufficiently large constant\\
	&&
\end{tabular}
\end{center}
\end{adjustwidth}

The following are decision variables in the optimization models:
\begin{adjustwidth}{0.25in}{0.25in} 
\begin{center}
\begin{tabular}{l c p{4.2in}}
	&&\\
	$y_i$ & : & spend of activity $i$\\
	$x_i$ & : & spend change from the baseline level of activity $i$\\
	$z_i$ & : & 1 if and only if $x_i$ is within region $S$ (i.e., $x_i$ equals to 0); 0, otherwise\\
	$z_i^L$ & : & 1 if and only if $x_i$ is within region $L$; 0, otherwise\\
	$z_i^R$ & : & 1 if and only if $x_i$ is within region $R$; 0, otherwise\\
	$z_i^{LR}$ & : & 1 if and only if $x_i$ is within regions $L$ or $R$; 0, otherwise\\
	&&
\end{tabular}
\end{center}
\end{adjustwidth}

To formulate our marketing resource allocation problem with the minimum spend change and maximum number of activities with spend change requirements, we first model it as a mixed integer nonlinear programming problem as follows.
\begin{subequations} \label{opt_MINLP}
\begin{align}
	\text{(MINLP)} \quad \max_{x, z} \quad & \sum_{i \in I} \left( \theta_i x_i^2 + \phi_i x_i + \psi_i \right) \label{opt_MINLP_obj} \\
	\textrm{s.t.} \quad & l_i - s_i \leq x_i \leq u_i - s_i, \: \forall i \in I \label{opt_MINLP_con1} \\
	& \sum_{i \in I} x_i \leq \left( \rho - 1 \right) \sum_{i \in I} s_i \label{opt_MINLP_con2} \\
	& \delta_i z_i \leq | x_i | \leq M z_i, \: \forall i \in I \label{opt_MINLP_con3} \\
	& \sum_{i \in I}{z_i} \leq m \label{opt_MINLP_con4} \\
	& z_i \in B, \: \forall i \in I \label{opt_MINLP_con5}
\end{align}
\end{subequations}
In the above formulation, objective function~(\ref{opt_MINLP_obj}) maximizes the total revenue generated from marketing activities. Constraint~(\ref{opt_MINLP_con1}) ensures the spend is within the bounds for each activity; constraint~(\ref{opt_MINLP_con2}) assures the total spend is within the upper bound of total budget; constraint~(\ref{opt_MINLP_con3}) stipulates the minimum change in spend; constraint~(\ref{opt_MINLP_con4}) guarantees the maximum number of changes; and, finally, constraint~(\ref{opt_MINLP_con5}) defines the integrality of decision variables.

The absolute value operator in constraint~(\ref{opt_MINLP_con3}) breaks the continuous range of $[l_i - s_i, u_i - s_i]$ for $x_i$ into at most two separate ranges and one point. To handle such discontinuity, we introduce two more binary variables and pre-calculate the bounds of the separate ranges to remove the absolute value operator from the constraint. Consequently, we reformulate the optimization problem as a mixed integer quadratic programming as follows.
\begin{subequations} \label{opt_MIQP}
\begin{align}
	\text{(MIQP)} \quad \max_{x, z^L, z^R} \quad & \sum_{i \in I} \left( \theta_i {x_i}^2 + \phi_i x_i + \psi_i \right) \label{opt_MIQP_obj}\\
	\textrm{s.t.} \quad & \sum_{i \in I} x_i \leq \left( \rho - 1 \right) \sum_{i \in I} s_i \label{opt_MIQP_con1}\\
	& l_i^L z_i^L + l_i^R z_i^R \leq x_i \leq u_i^L z_i^L + u_i^R z_i^R, \: \forall i \in I \label{opt_MIQP_con2}\\
	& z_i^L + z_i^R \leq 1, \: \forall i \in I \label{opt_MIQP_con3}\\
	& \sum_{i \in I}{\left( z_i^L + z_i^R \right)} \leq m \label{opt_MIQP_con4}\\
	& z_i^L \in B_i^L, z_i^R \in B_i^R, \: \forall i \in I \label{opt_MIQP_con5}
\end{align}
\end{subequations}
Compared to the MINLP, constraints~(\ref{opt_MINLP_con1}), (\ref{opt_MINLP_con3})-(\ref{opt_MINLP_con5}) are replaced by constraints~(\ref{opt_MIQP_con2})-(\ref{opt_MIQP_con5}) in the MIQP. Note that we pre-calculate the bounds of regions $L$ and $R$ in constraint~(\ref{opt_MIQP_con2}) and restrict the values of $z_i^L$ and $z_i^R$ in constraint~(\ref{opt_MIQP_con5}) based on the relationships established in Table~\ref{tab:MC_change}. Also, constraint (\ref{opt_MINLP_con1}) is no longer needed as it is absorbed into constraint~(\ref{opt_MIQP_con2}).

An alternative to handling the requirement of minimum change in spend (other than explicitly considering all possible feasible regions for an activity spend) is to directly restrict the square of the spend change to be no less than the square of the minimum change (i.e., ${x_i}^2 \geq {\delta_i}^2$) if there is a change. However, such a quadratic constraint has a convex function and ``$\geq$'' sign, resulting in a non-concave maximization optimization problem~\cite{Boyd_et_al_2004}.

Applying the perspective reformulation to an optimization problem containing semi-continuous variables with multiple continuous ranges may improve the solution process. In addition to the two binary variables $z_i^L$ and $z_i^R$, we introduce one more binary variable $z_i^{LR}$ to indicate if the activity spend changes from its baseline level (i.e., $y_i \neq s_i$ or $x_i \neq 0$). By replacing the objective function with its perspective function (i.e., $f_i (x_i) \to z_i^{LR} f_i (x_i / z_i^{LR})$) and adding one more constraint (i.e., $z_i^{L} + z_i^{R} = z_i^{LR}$) in model (\ref{opt_MIQP}), we obtain the perspective reformulation.

However, it is well known that applying the perspective reformulation may tighten the continuous relaxation at the expense of increasing the nonlinearity of the objective function. By fully exploiting the property that our objective function is a (concave) separable quadratic function, we further reformulate the model as a mixed integer second-order cone program. Considering the (concave) separable quadratic function, we have:
\begin{align}
	z_i^{LR} f_i \left( \frac{x_i}{z_i^{LR}} \right) & = \theta_i \left( \frac{{x_i}^2}{z_i^{LR}} \right) + \phi_i x_i + \psi_i. \nonumber
\end{align}
We can linearize $\theta_i ({x_i}^2 / z_i^{LR})$ by replacing it with a new variable times a minus sign ($- e_i$) and adding one more quadratic constraint ($- e_i z_i^{LR}\leq \theta_i {x_i}^2$) so that the model becomes a mixed integer second-order cone program. Note that when $z_i^{LR} = 0$, $x_i = 0$ and $\theta_i ({x_i}^2 / z_i^{LR}) = 0$, and the new quadratic constraint is trivially satisfied ($- e_i 0 \leq \theta_i 0$) and leaves $e_i$ unconstrained. However, since we are maximizing a concave objective function, we need to add a non-negativity constraint to the new variable ($e_i \geq 0$) to avoid the objective function being unbounded. The reformulated model is as follows.
\begin{subequations} \label{opt_MISOCP}
\begin{align}
	\textrm{(MISOCP)} \quad \max_{x, e, z^L, z^R, z^{LR}} \quad & \sum_{i \in I} \left( - e_i + \phi_i x_i + \psi_i \right) \label{opt_MISOCP_obj}\\
	\textrm{s.t.} \quad & - e_i z_i^{LR}\leq \theta_i {x_i}^2, \: \forall i \in I \label{opt_MISOCP_con1}\\
	& \sum_{i \in I} x_i \leq \left( \rho - 1 \right) \sum_{i \in I} s_i \label{opt_MISOCP_con2}\\
	& l_i^L z_i^L + l_i^R z_i^R \leq x_i \leq u_i^L z_i^L + u_i^R z_i^R, \: \forall i \in I \label{opt_MISOCP_con3}\\
	& z_i^L + z_i^R \leq 1, \: \forall i \in I \label{opt_MISOCP_con4}\\
	& \sum_{i \in I} \left( z_i^L + z_i^R \right) \leq m \label{opt_MISOCP_con5}\\
	& z_i^L + z_i^R = z_i^{LR}, \: \forall i \in I \label{opt_MISOCP_con6}\\
	& z_i^L \in B_i^L, z_i^R \in B_i^R, \: \forall i \in I \label{opt_MISOCP_con7}\\
	& e_i \geq 0, \: \forall i \in I \label{opt_MISOCP_con8}
\end{align}
\end{subequations}

\section{Numerical Experiment and Real-World Application}
We conduct a series of computational tests, including a set of numerical instances and a real-world application, to compare the solution process of the MIQP and MISOCP formulations using Gurobi solver.

\subsection{Numerical Instances}
Although there are several existing libraries of quadratic programming instances (e.g., \cite{Maros_et_al_1999, Furini_et_al_2019}) available, the structures of such instances are different from the one discussed in this paper. Hence, we construct instances that are more consistent with our model~(\ref{opt_MIQP}). Specifically, we limit the total activity budget to be no more than 1\% away from the total baseline activity budget (i.e., $\rho = 1.01$), add two more linear constraints to the model (i.e., $\sum_{i \in I}{\tau_i x_i} \leq (\zeta - 1) \sum_{i \in I}{\tau_i s_i}$ and $\sum_{i \in I}{\gamma_i x_i} \geq (\eta - 1) \sum_{i \in I}{\gamma_i s_i}$), and restrict the spend change for each activity to be no less than $\epsilon$ times the current spend (i.e., $\delta_i = \epsilon s_i$).

In order to generate instances at different difficulty levels, we employ the method used in~\cite{Pisinger_2005}, which investigated the hard instances of the Knapsack Problem (KP), to generate three sets of objective function coefficients based on the two added linear constraints: (i) uncorrelated, (ii) weakly correlated, and (iii) strongly correlated. We imitate these two constraints as two knapsack constraints in spite of that the other decision (i.e., change amount) also needs to be made if an activity is chosen for spend change. The relationship between the objective function coefficients ($\theta_i$, $\phi_i$ and $\psi_i$) and constraint coefficients ($\tau_i$ and $\gamma_i$) of these three sets are as suggested per~\cite{Pisinger_2005} as follows.

\begin{itemize}
\item \sloppy (i) \textit{Uncorrelated}: $\theta_i$ is uniformly distributed in $[-10, -1]$, and each of $\phi_i$ and $\psi_i$ is uniformly distributed in $[1, 10]$. There is no correlation between the objective function coefficients ($\theta_i, \phi_i, \psi_i$) and the constraint coefficients ($\tau_i, \gamma_i$). Such instances are supposedly relatively easy to solve from the perspective of KP since its ratio of profit to weight has quite a wide variation.

\item \sloppy (ii) \textit{Weakly correlated}: given $\tau_i$ and $\gamma_i$, $\theta_i$ is uniformly distributed in $\left [ - \left( \frac{\tau_i + \gamma_i}{2} + 1 \right), - \left( \frac{\tau_i + \gamma_i}{2} - 1 \right) \right ]$ and each of $\phi_i$ and $\psi_i$ is distributed in $\left[ \left( \frac{\tau_i + \gamma_i}{2} - 1 \right), \left( \frac{\tau_i + \gamma_i}{2} + 1 \right) \right]$. There is a strong correlation between each of the objective function coefficients $\theta_i$, $\phi_i$ and $\psi_i$ and the constraint coefficients $\tau_i$ and $\gamma_i$. Such instances are supposedly harder to solve since its ratio of profit to weight has much smaller variation.

\item \sloppy (iii) \textit{Strongly correlated}: given $\tau_i$ and $\gamma_i$, $\theta_i$ is computed as $- \left( \frac{\tau_i + \gamma_i}{2} + 1 \right)$ and each of $\phi_i$ and $\psi_i$ is computed as $\left( \frac{\tau_i + \gamma_i}{2} + 1 \right)$. Each of the objective function coefficients $\theta_i$, $\phi_i$ and $\psi_i$ is a function of the constraint coefficients $\tau_i$ and $\gamma_i$. Such instances are supposedly the hardest among the three sets to solve since its ratio of profit to weight has no variation.
\end{itemize}

The steps to generate data for the numerical instances are: (1) randomly choose $l_i$ and $u_i$ from uniform distribution of $[1, 5]$ and $[5, 10]$, respectively; (2) randomly choose $s_i$ from a uniform distribution of $[l_i, u_i]$; (3) randomly choose $\tau_i$ from a uniform distribution of $[1, 10]$, same as to $\gamma_i$; (4) generate the three sets of uncorrelated/correlated $\theta_i$, $\phi_i$ and $\psi_i$ using $\tau_i$ and $\gamma_i$ accordingly; (5) randomly choose $\zeta$ and $\eta$ from uniform distribution of $[0.90, 1.00]$ and $[1.00, 1.10]$, respectively; (6) set $\rho = 1.01$; and, (7) compute $l_i^L$, $u_i^L$, $l_i^R$ and $u_i^R$ using $s_i$ and $\delta_i$ according to Table~\ref{tab:MC_change}. All the generated random variates are rounded to two decimal places.

In addition to the three sets of uncorrelated/correlated objective function coefficients, we parameterize the total number of activities ($n = \{ 500, 750, 1000 \}$), the minimum change in spend (in percentage of current spend, i.e., $\epsilon = \{ 5\%, 10\%, 20\% \}, \delta_i = \epsilon s_i$), and the maximum number of changes (in percentage of total activities, i.e., $\xi = \{ 50\%, 75\%, 100\% \}, m = \xi n$). For each combination of correlation, $n$, $\epsilon$, and $\xi$, we generate 35 instances, resulting in 2,835 instances in total. The instances are solved using Amazon Elastic Compute Cloud (EC2) (running Amazon Linux AMI 2018.03 with dual Intel(R) Xeon(R) 16-core E5-2686 v4 CPU @ 2.3GHz and 244GB memory) with Gurobi version 7.0.2. Moreover, we set the following Gurobi configuration: the maximum run time is set to 100 seconds, the acceptable optimality gap is set to 0, the computing thread is set to 4, and other parameters are at their default values. For performance comparison between the MIQP and MISOCP formulations, we report the median of optimality gap, ratio of Gurobi run time and ratio of number of nodes from the branch-and-bound scheme employed by Gurobi. We chose median over mean to avoid the disproportionate influence of outliers.

\subsection{Numerical Results}
We report the numerical experiment results in Tables~\ref{tab:NumericalExperimentPerformanceComparison} as follows, which shows the problem characteristics for the MIQP and MISOCP formulations grouped by all, the correlation, the total number of activities ($n$), the minimum change in activity spend in percentage of the current spend ($\epsilon$), and the maximum number of spend changes in percentage of the total activities ($\xi$).
\begin{table}[ht!]
\footnotesize
\caption{Performance comparison of the MIQP and MISOCP formulations for the numerical experiment. We group the result by all, correlation, total number of activities ($n$), minimum change in activity spend (in percentage of current spend, i.e., $\delta_i = \epsilon s_i$) and maximum number of spend changes (in percentage of total activities, i.e., $m = \xi n$), and report the median to avoid outlier influence. Both run time and node count ratios use MIQP result as the base, e.g., overall the MISOCP spends 91\% less time and requires 99\% fewer nodes than the MIQP.}
\setlength{\extrarowheight}{0.1in}
\begin{adjustwidth}{-1in}{-1in} 
\begin{center}
	\begin{tabular}{ c | c | c  c | c  c }
		\hlineB{4}
		\multirow{2}{*}{Group} & \multirow{2}{*}{Instances} & MIQP & MISOCP & Run time & Node count \\
		&  & Optimality gap & Optimality gap & ratio & ratio \\
		\hlineB{4}
		All & 2,835 & 4.08\% & 0.00\% & 0.09 & 0.01 \\
		\hline
		Uncorrelated & 945 & 3.96\% & 0.00\% & 0.09 & 0.01 \\
		Weakly correlated & 945 & 4.16\% & 0.00\% & 0.08 & 0.01 \\
		Strongly correlated & 945 & 4.24\% & 0.00\% & 0.09 & 0.01 \\
		\hline
		$n$ = 500 & 945 & 2.74\% & 0.00\% & 0.04 & 0.00 \\
		$n$ = 750 & 945 & 4.32\% & 0.00\% & 0.09 & 0.01 \\
		$n$ = 1,000 & 945 & 5.22\% & 0.00\% & 0.14 & 0.03 \\
		\hline
		$\epsilon$ = 5\% & 945 & 0.96\% & 0.00\% & 0.16 & 0.02 \\
		$\epsilon$ = 10\% & 945 & 2.47\% & 0.00\% & 0.10 & 0.01 \\
		$\epsilon$ = 20\% & 945 & 9.12\% & 0.00\% & 0.05 & 0.01 \\
		\hline
		$\xi$ = 50\% & 945 & 10.03\% & 0.00\% & 0.11 & 0.02 \\
		$\xi$ = 75\% & 945 & 2.32\% & 0.00\% & 0.09 & 0.01 \\
		$\xi$ = 100\% & 945 & 1.70\% & 0.00\% & 0.08 & 0.01 \\
		\hlineB{4}
	\end{tabular}
\end{center}
\end{adjustwidth}
\label{tab:NumericalExperimentPerformanceComparison}
\end{table}

Out of the 2,835 instances, 1,995 (or 70\%) of the MISOCP can be solved to optimality within 100 seconds, compared to 37 (or 1\%) of the MIQP. Furthermore, the median optimality gap for the MIQP overall is 4.08\%. The gap tends to grow with the increase of correlation, $n$ and $\epsilon$, and the decrease of $\xi$, among which the gap varies most (from 10.03\% to 1.70\%) with the change in $\xi$ (from 50\% to 100\%). In addition, in general the MISOCP spends 91\% less time and requires 99\% fewer nodes than the MIQP to obtain a solution with a much smaller optimality gap. In particular, the MISOCP spends 89\% less time and requires 98\% fewer nodes than the MIQP when $\xi = 50\%$, which represents the hardest instances (with the largest optimality gap of 10.03\%) in our experiment. Note that the run time, optimality gap and node count together may better reflect the difficulty of the instances than any one of them alone.

\subsection{A Real-World Application}
We discuss our experience of applying the MISOCP method to an instance of adjusting the marketing plan at one of the largest independent food retailers in the United States. The company employs marketing instruments that can be categorized into ten different groups, specified as TV ads, direct mail, check-out coupon, digital media, digital coupon, flyer, outdoor ads, radio ads, paper coupon, and in-store displays. Among these groups, some are further divided into several subgroups, each reflecting an advertising campaign targeted at a local market. For example, a TV ad in a back-to-school campaign targeting a specific city station is considered as part of a subgroup. There are close to 200 marketing instrument groups/subgroups considered, and each of them may produce a sales lift on certain products or product categories in certain store groups for some time periods. The numbers of involved store groups and product categories are both in the hundreds.

A specific marketing instrument subgroup on a group of stores in a specific week is considered a marketing activity, and the resources (marketing spend) allocated to specific marketing activities constitute decisions to be made in a marketing mix optimization based planning process. The product category sales lifts by marketing instrument, store group, and week combinations are given by separable quadratic functions that are estimated based on historical data. The objective for the marketing mix optimization is to maximize the total estimated sales by optimally allocating the total marketing budget to marketing activities, subject to certain business rules. Typical business rules include but are not limited to: (1) lower and upper bounds on the spend of each marketing activity for each week; (2) lower and upper bounds on the total spend of each pair of marketing instrument-store groups for a part or the entire planning horizon; (3) lower and upper bounds on the total spend of a specific combination of marketing activities (e.g., total TV spend including all TV subgroups across all weeks, total TV and radio spends across all weeks, etc.); and (4) lower and upper bounds on the total spend of all marketing activities. The baseline marketing mix optimization model involves about 3,000 valid marketing activities (decision variables) and 2,600 business rules (constraints).

After the baseline marketing mix optimization problem is solved as a part of the company's annual planning process, an adjustment of the spend on each marketing activity is then contemplated as part of periodic marketing plan reviews. In one scenario of this adjustment optimization, we allow the maximum number of activities with spend change to be 25\% of the total marketing activities, set the minimum change in each spend to be 10\% of its baseline spend, and limit the total budget to be 90\% of its baseline level. In addition, we relax the budgets of collective marketing activity spends stipulated in other business rules to be from 80\% to 120\% of their baseline levels. The adjustment optimization problem aims to maximize the total estimated sales subject to the new spend adjustment requirements.

We apply both MIQP and MISOCP formulations to solve the aforementioned adjustment optimization problem. Compared to our baseline planning model, which consists of nearly 3,000 decision variables and 2,600 constraints, the spend adjustment model adds additional binary decision variables (about 6,000 to the MIQP and 9,000 to MISOCP) and constraints (about 12,000 to the MIQP and 18,000 to MISOCP) to enforce the requirements of minimum change in spend and maximum number of changes. By using Gurobi on the same computer mentioned in Section~4.1, with the same configuration except for extending the maximum run time to 600 seconds, we obtain the results in Table~\ref{tab:RealWorldApplicationPerformanceComparison}.
\begin{table}[ht!]
\footnotesize
\caption{Performance comparison of the MIQP and MISOCP formulations for the real-world application, adjusting activity budgets for a marketing mix optimization. The problem consists of nearly 9,000 marketing activities and 2,600 business rule constraints.}
\setlength{\extrarowheight}{0.1in}
\begin{adjustwidth}{-1in}{-1in} 
\begin{center}
	\begin{tabular}{ c  c | c  c | c  c }
		\hlineB{4}
		\multicolumn{2}{ c |}{Optimality gap} & \multicolumn{2}{ c |}{Run time (seconds)} & \multicolumn{2}{c}{Node count} \\
		MIQP & MISOCP & MIQP & MISOCP & MIQP & MISOCP \\
		\hlineB{4}
		0.0120\% & 0.0000\% & 600.00 & 1.38 & 3.4e5 & 0 \\
		\hlineB{4}
	\end{tabular}
\end{center}
\end{adjustwidth}
\label{tab:RealWorldApplicationPerformanceComparison}
\end{table}

\begin{figure}[h!]
  \includegraphics[width=\linewidth]{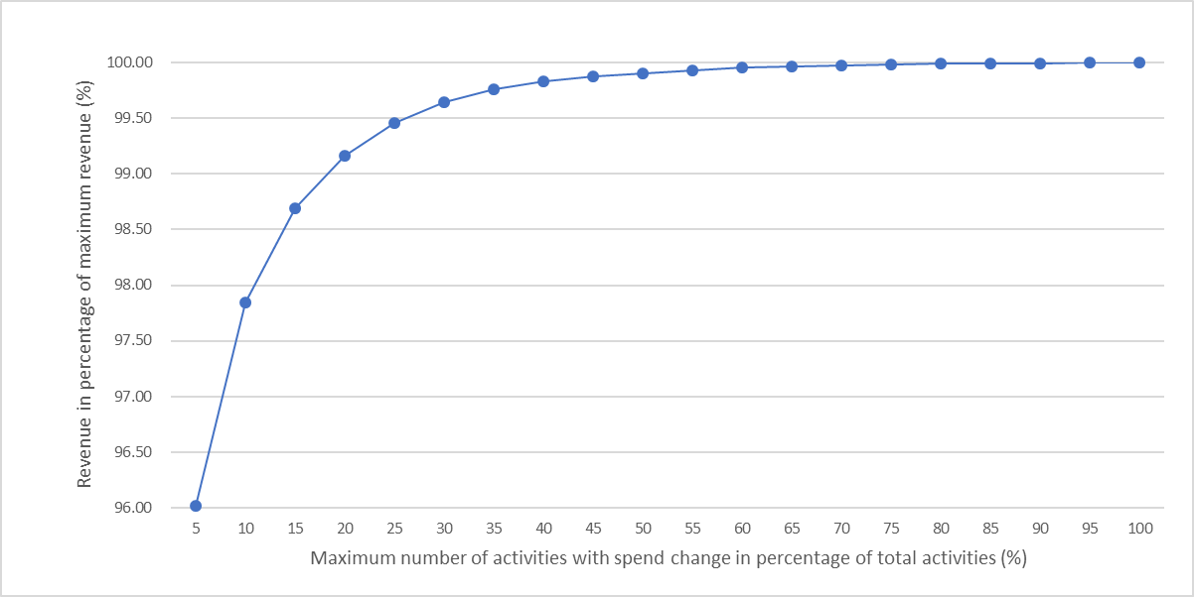}
  \caption{Pareto curve of revenue (in percentage of maximum revenue) against maximum number of activities with spend change (in percentage of total number of activities) for the real-world application.}
  \label{fig:RealWorldApplicationParetoCurve}
\end{figure}

In this real-world application, the MIQP provides a good heuristic solution, which is almost as good as the true optimal solution, in a little more of ten seconds and spends the rest of the time shrinking the optimality gap to 0.0120\%. On the other hand, the MISOCP obtains the optimal solution within two seconds without even exploring the branch-and-bound tree. With the rapid computation time resulting from the use of  MISOCP, we were able to produce a Pareto curve of the revenue (in percentage of maximum revenue, which is obtained by allowing all marketing activities subject to spend change) against the maximum number of activities allowed to change (in percentage of total activities). We illustrate such a Pareto curve in Figure~\ref{fig:RealWorldApplicationParetoCurve}. The curve provides insights regarding the number of spend adjustment versus its revenue gain. For example, it clearly shows that the spend adjustment of 25\% of total activities is around its knee point, translating into 99.5\% of maximum revenue. Management can utilize the Pareto curve to explore the tradeoffs between the activity adjustments and revenue opportunities by further considering available resources, costs of adjustment and desired benefits. In addition to swiftly analyzing various scenarios, the MISOCP also enables us to expand our marketing mix adjustment optimization problem at a more granular level, indicating more achievable benefits in practice.

\section{Discussion}
The outcomes from our numerical experiment show that the MISOCP formulation significantly outperforms the MIQP formulation from the perspective of optimality gap, run time and node count. For the majority of the instances of various problem sizes and difficulty levels we tested, the MISOCP formulation can be solved to optimality in tens of seconds while the MIQP formulation only resulted in near-optimal solutions with a much longer running time. Given the maximum run time of 100 seconds, for very few instances (less than 5.5\% of total) the MIQP do outperform the MISOCP, by either providing a smaller optimality gap or spending less solution time. Although for these instances the MIQP branch-and-bound tree sizes were much larger than the corresponding MISOCP one, the MIQP formulation was solved utilizing re-optimization capabilities of the branch-and-bound algorithm to shorten the overall solution time.

In general, the MISOCP formulation provides a much tighter continuous relaxation compared to the MIQP formulation at the expense of a slightly larger problem size (i.e., more decision variables and constraints). Due to the tighter continuous relaxation, the MISOCP formulation requires much smaller brand-and-bound tree sizes to reach a given optimality gap. If we specify an acceptable optimality gap of 5\%, in addition to the allowed run time of 100 seconds in the numerical experiment, the MISOCP has, in median, a 66\% smaller actual optimality gap, a 83\% shorter run time, and a 100\% smaller resultant tree size as compared to the MIQP formulation. Moreover, the MIQP formulation finds a heuristic solution of a single-digit optimality gap quite fast but has much difficulty in further reducing the gap. If we specify different allowed run times in the numerical experiment, the MIQP reaches a median optimality gap of 5.49\% in 10 seconds, 4.43\% in 50 seconds and 4.08\% in 100 seconds. For many instances, it would take the MIQP hours, if not days, to reach an optimality gap of 1\%. In fact, the MIQP swiftly provides a heuristic solution which is in proximity of the optimal, sometimes even the optimal itself, with a noticeable optimality gap, and spends the rest of the time proving the solution quality.

When using Gurobi with default configuration to solve the problem in both formulations, even though the MIQP formulation had the advantage of re-optimization capabilities of the branch-and-bound algorithm, the MISOCP formulation still resulted in better performance. In addition to a much smaller tree size, the outstanding performance of the MISOCP is also attributed to an efficient algorithm Gurobi employed using a linearized, outer approximation approach~\cite{Gurobi_MIQCP_2020}, which may re-use information from previously solved nodes in its branch-and-bound tree. If we solve a continuous quadratic constraint programming (QCP) at each node (by configuring Gurobi appropriately) when solving the MISOCP~\cite{Gurobi_Document_2020}, then the solution time would be much longer (even longer than the MIQP).

The computing speed gained with the MISOCP formulation directly translates into practical implications as evidenced by the real-world application in Section~4.3. While the business rule on the maximum number of changes may be dictated by available resources, contractual obligations and other factors, using the Pareto between revenue improvements and number of activities with marketing spend change, we are able to provide management with a data-driven approach to determine the maximum number of activities to allow for spend changes. In fact, with some flexibility and information on costs of implementing marketing spend changes, we should be able to recommend an optimal number of activities with marketing spend change that will generate the best marketing returns.

\section{Conclusions}
In this paper, we address an important variant of the marketing mix optimization problem that adjusts budget allocations for a large number of marketing activities. For practical reasons, the optimization problem needs to consider the minimum change in activity spend as well as the maximum number of activities with spend change. We provide a systematic way to model the minimum change requirement. Moreover, we develop a method to apply perspective reformulation on a semi-continuous variable with absolute value operator, and further on a variable with several separate continuous ranges as its feasible regions. This method aims to obtain a tighter continuous relaxation. By using this relaxation technique and exploiting the special property of the (concave) separable quadratic objective function, we reformulate the original optimization problem as a mixed integer second-order cone programming (MISOCP) to accelerate the speed of solution process. Based on the computational tests involving numerous specially designed instances and the real-world application, the proposed reformulation improves significantly by requiring less computation times and/or providing better solutions. The improvement can be attributed to the smaller branch-and-bound tree and much efficient Gurobi algorithm for MISOCP. The technique and reformulation discussed in this paper can be readily applied to other problems involving semi-continuous variables with absolute value operators and cardinality constraint, e.g., a mean-variance portfolio selection problem which allows both long selling and short selling with minimum buy-in threshold.

Our future research work will extend our approach to similar problems with non-separable quadratic objective functions. Note that, in some cases, the spend of a marketing activity may have a significant influence on the performance of other activities, resulting in a non-separable quadratic objective function. Some researchers have proposed several methods to decompose a non-separable Q matrix into two matrices (generic and diagonal) so that the perspective reformulation can be applied to semi-continuous variables for a tighter continuous relaxation. Our approach may be combined with the matrix decomposition methods to improve the solution process for problems with non-separable objective functions. This may provide a promising research direction to explore in the future.


\end{document}